\documentclass[10pt,reqno,twoside]{amsart}
\usepackage{amssymb,amsmath,amsthm,soul,color,paralist}
\usepackage{t1enc}
\usepackage{comment}
\usepackage[cp1250]{inputenc}
\usepackage{a4,indentfirst,latexsym}
\usepackage{graphics}
\usepackage{mathrsfs}
\usepackage{cite,enumitem,graphicx}
\usepackage[colorlinks=true,urlcolor=blue,
citecolor=red,linkcolor=blue,linktocpage,pdfpagelabels,
bookmarksnumbered,bookmarksopen]{hyperref}
\usepackage[english]{babel}
\usepackage[left=2.50cm,right=2.50cm,top=2.72cm,bottom=2.72cm]{geometry}
\usepackage[metapost]{mfpic}
\usepackage[colorinlistoftodos]{todonotes}
\usepackage[normalem]{ulem}

\numberwithin{equation}{section}

\allowdisplaybreaks

\def\Xint#1{\mathchoice
  {\XXint\displaystyle\textstyle{#1}}%
  {\XXint\textstyle\scriptstyle{#1}}%
  {\XXint\scriptstyle\scriptscriptstyle{#1}}%
  {\XXint\scriptscriptstyle\scriptscriptstyle{#1}}%
  \!\int}
\def\XXint#1#2#3{{\setbox0=\hbox{$#1{#2#3}{\int}$}
  \vcenter{\hbox{$#2#3$}}\kern-.5\wd0}}
\def\-int{\Xint -}

\newcommand{\R}{\mathbb{R}}

\newcommand{\p}{p^{*}}

\newcommand{\ri}{\rightarrow}

\DeclareMathOperator{\dive}{div}

\DeclareMathOperator{\e}{\varepsilon}

\newtheorem{lem}{Lemma}[section]
\newtheorem{thm}{Theorem}[section]

\newtheorem{remark}{Remark}[section]

\begin{document}
\title[the Choquard equation involving the $p$-Laplacian operator]
{Regularity and Pohozaev identity for the Choquard equation involving the $p$-Laplacian operator}

\author[V. Ambrosio]{Vincenzo Ambrosio}
\address{Vincenzo Ambrosio \hfill\break\indent
Dipartimento di Ingegneria Industriale e Scienze Matematiche \hfill\break\indent
Universit\`a Politecnica delle Marche\hfill\break\indent
Via Brecce Bianche, 12\hfill\break\indent
60131 Ancona (Italy)}
\email{v.ambrosio@univpm.it}





\maketitle

\begin{abstract}
In this paper, we study the regularity of weak solutions for a class of nonlinear Choquard equations driven by 
the $p$-Laplacian operator. We also establish a Pohozaev type identity. 
\end{abstract}

\section{Introduction}
Let $p\in [2, \infty)$ and $N>p$.
Let us consider the following nonlinear problem:
\begin{equation}\label{P}
\left\{
\begin{array}{ll}
-\Delta_{p} u+|u|^{p-2}u=(I_{\alpha}*F(u))f(u) \, \mbox{ in } \R^{N}, \\
u\in W^{1, p}(\R^{N}), 
\end{array}
\right.
\end{equation}
where $\Delta_{p}u:=\dive(|\nabla u|^{p-2}\nabla u)$ is the $p$-Laplacian operator, $\alpha\in (0, N)$, $I_{\alpha}:\R^{N}\setminus\{0\}\ri \R$ is the Riesz potential given by
\begin{align*}
I_{\alpha}(x)=\frac{\Gamma(\frac{N-\alpha}{2})}{\Gamma(\frac{\alpha}{2})\pi^{\frac{N}{2}}2^{\alpha}|x|^{N-\alpha}},
\end{align*}
the nonlinearity $f:\R\ri \R$ is a continuous function satisfying the following conditions:
\begin{compactenum}[$(f_1)$]
\item $\displaystyle{\lim_{|t|\ri 0}\frac{|f(t)|}{|t|^{\frac{N+\alpha}{N}\frac{p}{2}-1}}=0}$,
\item $\displaystyle{\limsup_{|t|\ri \infty}\frac{|f(t)|}{|t|^{\frac{N+\alpha}{N-p}\frac{p}{2}-1}}<\infty}$,
\end{compactenum}
and $F(t):=\int_{0}^{t}f(\tau)\, d\tau$.
We recall that when $p=2$, $N=3$ and $F(t)=\frac{t^{2}}{2}$ in \eqref{P}, then we obtain the following Choquard-Pekar equation
\begin{align}\label{CHOQ}
-\Delta u+u=(I_{2}*|u|^{2})u \, \mbox{ in } \R^{3},
\end{align}
which seems to originate from Fr\"ohlich and Pekar's model of the quantum polaron and then used by Choquard 
to describe an electron trapped in its own hole, in a certain approximation to Hartree-Fock theory of one component plasma.
Equation \eqref{CHOQ} also appears as a model of self-gravitating matter and it is usually called the nonlinear Schr\"odinger-Newton equation. 
For more details on Choquard type equations and some of its variants and generalizations, we refer to \cite{MVSG} and the references therein.

We stress that the choice of the exponents $\frac{N+\alpha}{N}\frac{p}{2}$ and $\frac{N+\alpha}{N-p}\frac{p}{2}$ in $(f_1)$ and $(f_2)$ is motivated by the following Hardy-Littlewood-Sobolev inequality:
\begin{thm}\label{HLS}\cite{LL}
Let $r, t\in (1, \infty)$ and $\mu\in (0, N)$ with $\frac{1}{r}+\frac{\mu}{N}+\frac{1}{t}=2$. Let $f\in L^{r}(\R^{N})$ and $h\in L^{t}(\R^{N})$. 
Then there exists a sharp constant $C(N, \mu, r)>0$, independent of $f$ and $h$, such that 
$$
\left|\iint_{\R^{2N}} \frac{f(x)h(y)}{|x-y|^{\mu}}\, dx dy\right|\leq C(N, \mu, r)\|f\|_{L^{r}(\R^{N})}\|h\|_{L^{t}(\R^{N})}.
$$
\end{thm}
Indeed, if $F(\tau)=|\tau|^{q}$ for some $q>0$, $\mu=N-\alpha$ and $u\in W^{1, p}(\R^{N})$, then Theorem \ref{HLS} implies that the term 
$$
\int_{\R^{N}} \left(\frac{1}{|x|^{\mu}}*F(u)\right) F(u)\, dx
$$
is well-defined whenever $F(u)\in L^{\tau}(\R^{N})$ with $\tau=\frac{2N}{N+\alpha}$. Since $W^{1, p}(\R^{N})$ is continuously embedded into $L^{r}(\R^{N})$ for all $r\in [p, \p]$, where $p^{*}:=\frac{Np}{N-p}$ is the critical Sobolev exponent, we have to require that $\tau q\in [p, \p]$ and so
$$
\frac{N+\alpha}{N}\frac{p}{2}\leq q\leq \frac{N+\alpha}{N-p}\frac{p}{2}.
$$
The endpoints in the above range of $q$ are extremal values for the Hardy-Littlewood-Sobolev inequality and sometimes called lower and upper Hardy-Littlewood-Sobolev critical exponents. 
We recall that if 
$\alpha\in (N-p, N)$ and $f$ satisfies $|f(t)|\leq C_{0}(|t|^{q_{1}-1}+|t|^{q_{2}-1})$ for all $t\in \R$, 
with $p< q_{1}\leq q_{2}<\frac{\alpha p}{N-p}$, 
then it is easy to verify that every weak solution to \eqref{P} belongs to $L^{\infty}(\R^{N})\cap C^{1, \sigma}_{loc}(\R^{N})$ for some $\sigma\in (0, 1)$ (see \cite[Theorem 2.2]{AY}).
In fact, the restrictions on $q_{1}$ and $q_{2}$ and Theorem \ref{HLS} permit to consider the convolution $K(x):=(I_{\alpha}*F(u))(x)$ as bounded term, and thus \cite[Theorem 1.11]{LI} can be applied for $-\Delta_{p}u+|u|^{p-2}u=K(x)f(u)$ in $\R^{N}$. In this paper, we prove that the same regularity holds for solutions of \eqref{P} by considering more general nonlinearities.
More precisely, the main result of this work is the following:
\begin{thm}\label{thm1}
Let $p\in [2, \infty)$, $N>p$ and $\alpha\in ((N-2p)_{+}, N)$. Assume that $f\in C(\R)$ fulfills $(f_1)$-$(f_2)$. 
Let $u\in W^{1, p}(\R^{N})$ be a weak solution to \eqref{P}, that is
\begin{align}\label{WF}
\int_{\R^{N}} |\nabla u|^{p-2}\nabla u \cdot \nabla \phi \, dx+\int_{\R^{N}} |u|^{p-2}u\phi \, dx=\int_{\R^{N}} (I_{\alpha}*F(u)) f(u)\phi\, dx
\end{align}
for all $\phi\in W^{1, p}(\R^{N})$. 
Then, $u\in L^{\infty}(\R^{N})\cap C^{1, \sigma}_{loc}(\R^{N})$ for some $\sigma\in (0, 1)$. Moreover, $u$ satisfies the following Pohozaev identity: 
\begin{align}\label{POHOZAEV}
\frac{N-p}{p} \int_{\R^{N}} |\nabla u|^{p}\, dx+\frac{N}{p} \int_{\R^{N}} |u|^{p}\, dx=\frac{N+\alpha}{2} \int_{\R^{N}} (I_{\alpha}*F(u))F(u) \, dx.
\end{align}
\end{thm}
To accomplish the boundedness of solutions to \eqref{P}, we use the growth conditions on $f$, Theorem \ref{HLS} and a Brezis-Kato argument \cite{BK} to reach the $L^{r}$-regularity for all $r\in [p, \infty)$, and finally we carry out a De Giorgi type iteration \cite{D}. The restriction $\alpha\in ((N-2p)_{+}, N)$ ensures that $\frac{N+\alpha}{N-p}\frac{p}{2}>p$, and that the term
$$ 
\left( \int_{\{|u|>\mu\}} |u|^{\p}\, dx \right)^{\frac{\alpha+2p-N}{2N}}
$$
goes to zero as $\mu\ri \infty$ and this allows us to have a convenient estimate to achieve the $L^{q}$-regularity; see formula \eqref{ROTHER}. 
The $C^{1, \sigma}_{loc}(\R^{N})$-regularity result is thus a direct consequence of the classical theorems established in \cite{DB, T}.
We recall that in the case $p=2$, Moroz and Van Schaftingen \cite{MVS} obtained a $W^{2, q}_{loc}(\R^{N})$-regularity result for weak solutions to \eqref{P} and involving a Berestycki-Lions type nonlinearity \cite{BL}, by proving a suitable nonlocal version of the Br\'ezis-Kato regularity result \cite[Theorem 2.3]{BK} and invoking the Calderon-Zygmund theory. The approach in \cite{MVS} does not seem to be so simple to adopt when the Laplacian operator $\Delta$ is replaced by the $p$-Laplacian operator $\Delta_{p}$ with $p\neq 2$, due to the nonlinear character of $\Delta_{p}$, and so we follow a different strategy. Finally, we show that every weak solution to \eqref{P} fulfills the Pohozaev identity \eqref{POHOZAEV}. For this purpose, we take advantage of a variational identity for locally Lipschitz continuous solutions of a general class of quasilinear equations demonstrated in \cite{DGMS}. 
Note that the identity \eqref{POHOZAEV} implies that \eqref{P} does not have nontrivial weak solutions in $W^{1, p}(\R^{N})$ when $F(t)=\frac{1}{q}|t|^{q}$ with $q\notin (\frac{N+\alpha}{N}\frac{p}{2}, \frac{N+\alpha}{N-p}\frac{p}{2})$. 
Moreover, \eqref{POHOZAEV} is in clear accordance with the case $p=2$ investigated in \cite{MVS}.

To our knowledge, this is the first time that the regularity of solutions to \eqref{P} and the Pohozaev identity \eqref{POHOZAEV} are considered in literature. Finally, we believe that Theorem \ref{thm1} will play a crucial role in studying existence, multiplicity and concentration of solutions for quasilinear Choquard equations with a general nonlinearity (see \cite{AAI}).

\section{Proof of Theorem \ref{thm1}}

Without loss of generality, we assume that $u\geq 0$ (otherwise, we can deal with $u^{+}:=\max\{u, 0\}$ and $u^{-}:=\max\{-u, 0\}$ separately).
Let $M>1$ and $u_{M}:=\min\{u, M\}$. Taking $u_{M}^{kp+1}$ ($k\geq 0$) as test function in \eqref{WF}, we get
\begin{align}\label{0AI}
\int_{\R^{N}} |\nabla u|^{p-2}\nabla u\cdot \nabla (u_{M}^{kp+1}) \, dx+\int_{\R^{N}} u^{p-1}u_{M}^{kp+1} \, dx=\int_{\R^{N}} (I_{\alpha}*F(u)) f(u)u_{M}^{kp+1}\, dx.
\end{align}
Now, using the Sobolev inequality, we see that
\begin{align}\begin{split}\label{1AI}
\int_{\R^{N}} |\nabla u|^{p-2}\nabla u\cdot \nabla (u_{M}^{kp+1}) \, dx&=\frac{kp+1}{(k+1)^{p}} \int_{\R^{N}} |\nabla (u_{M})^{k+1}|^{p}\, dx \\
&\geq S_{*} \frac{kp+1}{(k+1)^{p}} \left[\int_{\R^{N}} |u_{M}^{k+1}|^{\p}\, dx\right]^{\frac{p}{\p}},
\end{split}\end{align}
where $S_{*}=S_{*}(N, p)>0$ denotes the best Sobolev constant.
Clearly, since $u_{M}\leq u$, we have
\begin{align}\label{2AI}
\int_{\R^{N}} u^{p-1}u_{M}^{kp+1} \, dx\geq \int_{\R^{N}} |u_{M}^{k+1}|^{p}\, dx.
\end{align}
On the other hand, exploiting Theorem \ref{HLS} and that $(f_1)$-$(f_2)$ yield
\begin{align}\label{FG}
|F(t)|\leq C(|t|^{\frac{N+\alpha}{N} \frac{p}{2}}+|t|^{\frac{N+\alpha}{N-p} \frac{p}{2}}) \quad \mbox{ for all } t\in \R,
\end{align}
we obtain
\begin{align}\begin{split}\label{3AI}
\left|\int_{\R^{N}} (I_{\alpha}*F(u)) f(u)u_{M}^{kp+1}\, dx  \right|&\leq C\|F(u)\|_{L^{\frac{2N}{N+\alpha}}(\R^{N})} \|f(u)u_{M}^{kp+1}\|_{L^{\frac{2N}{N+\alpha}}(\R^{N})} \\
&\leq C\left(\|u\|^{p}_{L^{p}(\R^{N})}+\|u\|_{L^{\p}(\R^{N})}^{\p} \right)^{\frac{N+\alpha}{2N}} \|f(u)u_{M}^{kp+1}\|_{L^{\frac{2N}{N+\alpha}}(\R^{N})} \\
&=:\tilde{C} \, \|f(u)u_{M}^{kp+1}\|_{L^{\frac{2N}{N+\alpha}}(\R^{N})}.
\end{split}\end{align}
By $(f_1)$, fixed $\e>0$ there exists $\delta_{0}>0$ such that
\begin{align}\label{fN0}
|f(t)|\leq \e |t|^{\frac{N+\alpha}{N} \frac{p}{2}-1} \quad \mbox{ for all } 0<|t|<\delta_{0}.
\end{align}
Let $0<\delta<\min\{\delta_{0}, 1\}$. 
Using $(f_2)$, we can find $K_{0}>0$ and $C_{1}>0$ such that
\begin{align}\label{fNI}
|f(t)|\leq C_{1} |t|^{\frac{N+\alpha}{N-p} \frac{p}{2}-1} \quad \mbox{ for all } |t|>K_{0}.
\end{align}
Let $K>\max\{K_{0}, 1 \}$. 
Fix $\mu>K$. Since $f$ is continuous in $\R$, there exists $C_{\delta, \mu}>0$ such that
\begin{align}\label{W}
|f(t)|\leq C_{\delta, \mu} |t|^{\frac{N+\alpha}{N-p} \frac{p}{2}-1} \quad \mbox{ for all } \delta\leq |t|\leq \mu.
\end{align}
Note that, thanks to $u_{M}\leq u$, \eqref{fN0}, \eqref{fNI}, \eqref{W}, $(x+y+z)^{a}\leq x^{a}+y^{a}+z^{a}$ for all $x, y, z\geq 0$ and $a\in (0, 1)$, we find
\begin{align*}
&\left[\int_{\R^{N}} |f(u)u_{M}^{kp+1}|^{\frac{2N}{N+\alpha}}\, dx\right]^{\frac{N+\alpha}{2N}}\\
&=\left[\int_{\{|u|<\delta\}} |f(u)u_{M}^{kp+1}|^{\frac{2N}{N+\alpha}} \, dx+\int_{\{\delta<|u|<\mu\}} |f(u)u_{M}^{kp+1}|^{\frac{2N}{N+\alpha}} \, dx+\int_{\{|u|>\mu\}}  |f(u)u_{M}^{kp+1}|^{\frac{2N}{N+\alpha}}\, dx  \right]^{\frac{N+\alpha}{2N}} \\
&\leq \left[\e^{\frac{2N}{N+\alpha}}\int_{\{|u|<\delta\}} ||u|^{\frac{N+\alpha}{N} \frac{p}{2}+kp}|^{\frac{2N}{N+\alpha}} \, dx+C_{\delta, \mu}^{\frac{2N}{N+\alpha}} \int_{\{\delta<|u|<\mu\}} ||u|^{\frac{N+\alpha}{N-p}\frac{p}{2}+kp}|^{\frac{2N}{N+\alpha}} \, dx+C_{1}^{\frac{2N}{N+\alpha}} \int_{\{|u|>\mu\}} ||u|^{\frac{N+\alpha}{N-p}\frac{p}{2}+kp}|^{\frac{2N}{N+\alpha}} \, dx  \right]^{\frac{N+\alpha}{2N}} \\
&\leq  \e \left(\int_{\{|u|<\delta\}} |u|^{(\frac{N+\alpha}{N} \frac{p}{2}+kp)\frac{2N}{N+\alpha}} \, dx\right)^{\frac{N+\alpha}{2N}}+C_{\delta, \mu}  \left(\int_{\{\delta<|u|<\mu\}} |u|^{(\frac{N+\alpha}{N-p} \frac{p}{2}+kp)\frac{2N}{N+\alpha}} \, dx\right)^{\frac{N+\alpha}{2N}}\\
&\quad+C_{1} \left(\int_{\{|u|>\mu\}} |u|^{(\frac{N+\alpha}{N-p} \frac{p}{2}+kp)\frac{2N}{N+\alpha}} \, dx\right)^{\frac{N+\alpha}{2N}} \\
&\leq \e \left(\int_{\{|u|<\delta\}} |u|^{p} |u|^{\frac{2N}{N+\alpha}kp}  \, dx\right)^{\frac{N+\alpha}{2N}}+C_{\delta, \mu}  \left(\int_{\{|u|<\mu\}} |u|^{(\frac{N+\alpha}{N-p} \frac{p}{2}+kp)\frac{2N}{N+\alpha}} \, dx\right)^{\frac{N+\alpha}{2N}} \\
&\quad+C_{1} \left(\int_{\{|u|>\mu\}} |u|^{(\frac{N+\alpha}{N-p} \frac{p}{2}+kp)\frac{2N}{N+\alpha}} \, dx\right)^{\frac{N+\alpha}{2N}}.
\end{align*}
Since for $|u|<\delta<1$ we have $|u|^{\frac{2N}{N+\alpha} kp}\leq |u|^{kp}$ (note that $\frac{2N}{N+\alpha}>1$), we see that
\begin{align*}
\e \left(\int_{\{|u|<\delta\}} |u|^{p} |u|^{\frac{2N}{N+\alpha}kp}  \, dx\right)^{\frac{N+\alpha}{2N}}\leq \e \left(\int_{\{|u|<\delta\}} |u|^{(k+1)p}  \, dx\right)^{\frac{N+\alpha}{2N}}\leq \e \int_{\{|u|<\delta\}} |u|^{(k+1)p}  \, dx+\e,
\end{align*}
where we used $x^{a}\leq x+1$ for all $x\geq 0$ and $a\in (0, 1)$. Therefore, 
\begin{align}\begin{split}\label{4AI}
&\left[\int_{\R^{N}} |f(u)u_{M}^{kp+1}|^{\frac{2N}{N+\alpha}}\, dx\right]^{\frac{N+\alpha}{2N}}\\
&\leq \e \int_{\{|u|<\delta\}} |u|^{(k+1)p}  \, dx+\e+C_{\delta, \mu}  \left(\int_{\{|u|<\mu\}} |u|^{(\frac{N+\alpha}{N-p} \frac{p}{2}+kp)\frac{2N}{N+\alpha}} \, dx\right)^{\frac{N+\alpha}{2N}} \\
&\quad+C_{1} \left(\int_{\{|u|>\mu\}} |u|^{(\frac{N+\alpha}{N-p} \frac{p}{2}+kp)\frac{2N}{N+\alpha}} \, dx\right)^{\frac{N+\alpha}{2N}}. 
\end{split}\end{align}
Note that, because $0<\delta<1<M$, it holds $u_{M}=u$,  and so 
$$
\int_{\{|u|<\delta\}} |u|^{(k+1)p}  \, dx=\int_{\{|u|<\delta\}} |u_{M}|^{(k+1)p}  \, dx.
$$
Combining \eqref{0AI}, \eqref{1AI}, \eqref{2AI}, \eqref{3AI} and \eqref{4AI}, we obtain
\begin{align*}
& S_{*} \frac{kp+1}{(k+1)^{p}} \left[\int_{\R^{N}} |u_{M}|^{(k+1)\p}\, dx\right]^{\frac{p}{\p}}+(1-\tilde{C} \e) \int_{\{|u|<\delta\}} |u_{M}|^{(k+1)p}  \, dx\\
&\leq \tilde{C} \e+\tilde{C} C_{\delta, \mu} \left(\int_{\{|u|<\mu\}} |u|^{(\frac{N+\alpha}{N-p} \frac{p}{2}+kp)\frac{2N}{N+\alpha}} \, dx\right)^{\frac{N+\alpha}{2N}}+\tilde{C} C_{1} \left(\int_{\{|u|>\mu\}} |u|^{(\frac{N+\alpha}{N-p} \frac{p}{2}+kp)\frac{2N}{N+\alpha}} \, dx\right)^{\frac{N+\alpha}{2N}}.
\end{align*}
Letting first $\e\ri 0$ and then $M\ri \infty$, we arrive at
\begin{align}\begin{split}\label{8AI}
& S_{*} \frac{kp+1}{(k+1)^{p}} \left[\int_{\R^{N}} |u|^{(k+1)\p}\, dx\right]^{\frac{p}{\p}}\\
&\leq \tilde{C} C_{\delta, \mu} \left(\int_{\{|u|<\mu\}} |u|^{(\frac{N+\alpha}{N-p} \frac{p}{2}+kp)\frac{2N}{N+\alpha}} \, dx\right)^{\frac{N+\alpha}{2N}}+\tilde{C} C_{1} \left(\int_{\{|u|>\mu\}} |u|^{(\frac{N+\alpha}{N-p} \frac{p}{2}+kp)\frac{2N}{N+\alpha}} \, dx\right)^{\frac{N+\alpha}{2N}}.
\end{split}\end{align}
Now, we observe that
\begin{align}\begin{split}\label{9AI}
\left(\int_{\{|u|<\mu\}} |u|^{(\frac{N+\alpha}{N-p} \frac{p}{2}+kp)\frac{2N}{N+\alpha}} \, dx\right)^{\frac{N+\alpha}{2N}} 
&=\left(\int_{\{|u|<\mu\}} |u|^{(\frac{N+\alpha}{N-p} \frac{p}{2}-p+(k+1)p)\frac{2N}{N+\alpha}} \, dx\right)^{\frac{N+\alpha}{2N}} \\
&\leq \mu^{\frac{N+\alpha}{N-p} \frac{p}{2}-p}  \left(\int_{\{|u|<\mu\}} |u|^{(k+1)\frac{2Np}{N+\alpha}} \, dx\right)^{\frac{N+\alpha}{2N}} \\
&\leq \mu^{\frac{N+\alpha}{N-p} \frac{p}{2}-p}  \left(\int_{\R^{N}} |u|^{(k+1)\frac{2Np}{N+\alpha}} \, dx\right)^{\frac{N+\alpha}{2N}}.
\end{split}\end{align}
On the other hand, using H\"older's inequality,
\begin{align}\begin{split}\label{10AI}
\left(\int_{\{|u|>\mu\}} |u|^{(\frac{N+\alpha}{N-p} \frac{p}{2}+kp)\frac{2N}{N+\alpha}} \, dx\right)^{\frac{N+\alpha}{2N}} 
&=\left(\int_{\{|u|>\mu\}} |u|^{(\frac{N+\alpha}{N-p} \frac{p}{2}-p+(k+1)p)\frac{2N}{N+\alpha}} \, dx\right)^{\frac{N+\alpha}{2N}} \\
&= \left(\int_{\{|u|>\mu\}} |u|^{(\frac{N+\alpha}{N-p}\frac{p}{2}-p)\frac{2N}{N+\alpha}} |u|^{(k+1)\frac{2Np}{N+\alpha}} \, dx\right)^{\frac{N+\alpha}{2N}} \\
&\leq  \left( \int_{\{|u|>\mu\}} |u|^{\p}\, dx \right)^{\frac{\alpha+2p-N}{2N}} \left( \int_{\R^{N}} |u|^{(k+1)\p}\, dx \right)^{\frac{p}{\p}} \\
&=: D(\mu) \left( \int_{\R^{N}} |u|^{(k+1)\p}\, dx \right)^{\frac{p}{\p}}.
\end{split}\end{align}
In view of \eqref{8AI}, \eqref{9AI} and \eqref{10AI}, we deduce that
\begin{align*}
& S_{*} \frac{kp+1}{(k+1)^{p}} \left[\int_{\R^{N}} |u|^{(k+1)\p}\, dx\right]^{\frac{p}{\p}}\\
&\leq \tilde{C} C_{\delta, \mu} \, \mu^{\frac{N+\alpha}{N-p} \frac{p}{2}-p}  \left(\int_{\R^{N}} |u|^{(k+1)\frac{2Np}{N+\alpha}} \, dx\right)^{\frac{N+\alpha}{2N}}+\tilde{C} C_{1} D(\mu) \left( \int_{\R^{N}} |u|^{(k+1)\p}\, dx \right)^{\frac{p}{\p}}.
\end{align*}
Since $\alpha+2p-N>0$, $D(\mu)\ri 0$ as $\mu\ri \infty$, and we can choose $\mu>K$ sufficiently large such that
$$
0\leq D(\mu)<\theta \frac{S_{*}}{\tilde{C}C_{1}} \frac{kp+1}{(k+1)^{p}} \quad \mbox{ with } \theta\in (0, 1).
$$
Hence,
\begin{align*}
(1-\theta)S_{*} \frac{kp+1}{(k+1)^{p}} \left[\int_{\R^{N}} |u|^{(k+1)\p}\, dx\right]^{\frac{p}{\p}}\leq \tilde{C} C_{\delta, \mu} \, \mu^{\frac{N+\alpha}{N-p} \frac{p}{2}-p}  \left(\int_{\R^{N}} |u|^{(k+1)\frac{2Np}{N+\alpha}} \, dx\right)^{\frac{N+\alpha}{2N}},
\end{align*}
and so
\begin{align}\label{ROTHER}
\|u\|_{L^{(k+1)\p}(\R^{N})}\leq C_{*}^{\frac{1}{k+1}} \left[ \frac{k+1}{(kp+1)^{\frac{1}{p}}} \right]^{\frac{1}{k+1}} \|u\|_{L^{(k+1) \frac{2Np}{N+\alpha}}(\R^{N})},
\end{align}
for some $C_{*}>0$ depending on $k$. Set $q:=\frac{2Np}{N+\alpha}$ and observe that $\p>q$ thanks to $\alpha>N-2p$. Now a bootstrap argument can start: because $u\in L^{\p}(\R^N)$, we can apply \eqref{ROTHER} with $k+1=\frac{\p}{q}$ to see that $u\in L^{(k+1)\p}(\R^{N})=L^{\frac{(\p)^{2}}{q}}(\R^{N})$. We can then exploit again \eqref{ROTHER} and, after $m$ iterations, we obtain that $u\in L^{\p(\frac{\p}{q})^{m}}(\R^N)$ and thus $u\in L^{\nu}(\R^{N})$ for all $\nu\in [\p, \infty)$. Consequently, $u\in L^{\nu}(\R^{N})$ for all $\nu\in [p, \infty)$. Combining this fact with \eqref{FG}, and applying Young's inequality for convolutions, it is easy to verify that 
$K(x):=(I_{\alpha}*F(u))(x)\in C_{0}(\R^{N}):=\{v\in C(\R^{N}): |v(x)|\ri 0 \mbox{ as } |x|\ri \infty \}$. 
In particular, $K\in L^{\infty}(\R^{N})$. Then $u$ satisfies $-\Delta_{p}u=\psi(x, u):=-|u|^{p-2}u+K(x)f(u)$ in  $\R^{N}$,
with $|\psi(x, u)|\leq C(1+|u|^{r-1})$ for some $r\in (p, \p)$. Now we show that $u\in L^{\infty}(\R^N)$ by means of a De Giorgi type iteration. 
Let $\rho\geq \max\{1, \|u\|^{-1}_{L^{r}(\R^{N})}\}$ and $v:=(\rho \|u\|_{L^{r}(\R^{N})})^{-1}u$. Set $w_{k}:=(v-1+2^{-k})^{+}$ for $k\in \mathbb{N}$ and $w_{0}:=v^{+}$. Note that $w_{k}\in W^{1, p}(\R^{N})$ and $0\leq w_{k+1}\leq w_{k}$ a.e. in $\R^{N}$. Let $U_{k}:=\|w_{k}\|^{r}_{L^{r}(\R^{N})}$. Put $\Omega_{k}:=\{x\in \R^{N}: w_{k}(x)>0\}$. Observe that $\Omega_{k+1}\subset \{w_{k}>2^{-(k+1)}\}$, $v(x)<2^{k+1}w_{k}(x)$ for $x\in \Omega_{k+1}$ and $|\Omega_{k+1}|\leq 2^{(k+1)r}U_{k}$. Testing $-\Delta_{p}u=\psi(x, u)$ with $w_{k+1}$ and using the growth assumption on $\psi$, we can see that
\begin{align*}
\|\nabla w_{k+1}\|^{p}_{L^{p}(\R^{N})}&=\int_{\R^{N}} |\nabla v|^{p-2} \nabla v \cdot \nabla w_{k+1}\, dx=(\rho \|u\|_{L^{r}(\R^{N})})^{1-p} \int_{\R^{N}} |\nabla u|^{p-2} \nabla u \cdot \nabla w_{k+1}\, dx\\
&=(\rho \|u\|_{L^{r}(\R^{N})})^{1-p} \int_{\R^{N}} \psi(x, u) w_{k+1}\, dx \\
&\leq C(\rho \|u\|_{L^{r}(\R^{N})})^{1-p} \left[ \int_{\Omega_{k+1}\cap \{u\leq 1\}} w_{k+1}\, dx+\int_{\Omega_{k+1}\cap \{u>1\}} u^{r-1} w_{k}\, dx\right] \\
&\leq C(\rho \|u\|_{L^{r}(\R^{N})})^{1-p} \left[ 2^{-(k+1)}|\Omega_{k+1}|+(\rho \|u\|_{L^{r}(\R^{N})})^{r-1} 2^{(k+1)(r-1)} U_{k} \right] \\
&\leq 2C 2^{(k+1)(r-1)} (\rho \|u\|_{L^{r}(\R^{N})})^{r-p} U_{k},
\end{align*}
where we used $w_{k+1}(x)\leq 2^{-(k+1)}$ for $x\in \Omega_{k+1}\cap \{u\leq 1\}$ and $u(x)<(\rho \|u\|_{L^{r}(\R^{N})}) 2^{k+1}w_{k}(x)$ for $x\in \Omega_{k+1}$. In view of the above inequality and using the H\"older and Sobolev inequalities, we arrive at
\begin{align*}
U_{k+1}\leq \|w_{k+1}\|^{r}_{L^{\p}(\R^{N})} |\Omega_{k+1}|^{1-\frac{r}{\p}}\leq \bar{C}^{k} (\rho \|u\|_{L^{r}(\R^{N})})^{\frac{r^{2}}{p}-r} U_{k}^{1+\frac{r}{N}},
\end{align*}
for some $\bar{C}>1$ independent of $k$. Let $\eta:=\bar{C}^{-\frac{N}{r}}\in (0, 1)$. Choosing 
$$
\rho:=\max\{1, \|u\|^{-1}_{L^{r}(\R^{N})}, (\|u\|_{L^{r}(\R^{N})}^{\frac{r^{2}}{p}-r}\eta^{-1})^{\frac{1}{\gamma}}\},
$$ 
where $\gamma:=\frac{r^{2}}{N}+r-\frac{r^{2}}{p}>0$, by induction we can prove that $U_{k}\leq \frac{\eta^{k}}{\rho^{r}}$ for all $k\in \mathbb{N}\cup \{0\}$, from which $U_{k}\ri 0$ as $k\ri \infty$. Since $w_{k}\ri (v-1)^{+}$ a.e. in $\R^{N}$ as $k\ri \infty$ and $w_{k}\leq v\in L^{r}(\R^{N})$, it follows from the dominated convergence theorem that $\|(v-1)^{+}\|_{L^{r}(\R^{N})}=0$. As a result, $v\leq 1$ a.e. in $\R^{N}$, that is, $\|u\|_{L^{\infty}(\R^{N})}\leq \rho \|u\|_{L^{r}(\R^{N})}$.
Therefore, $u$ is a bounded solution to $-\Delta_{p}u=\psi=-|u|^{p-2}u+K(x)f(u)\in L^{\infty}(\R^{N})$, and according to \cite{DB, T} we can infer that $u\in C^{1, \sigma}_{loc}(\R^{N})$ for some $\sigma\in (0, 1)$.
\begin{remark}
If we assume that $u$ is a positive solution to \eqref{P}
and $f$ is an odd continuous function satisfying $f\geq 0$ in $(0, \infty)$, 
$f(t)=o(|t|^{p-1})$ as $|t|\ri 0$ instead of $(f_1)$, and $(f_2)$, then we can show that $u$ decays exponentially at infinity. To prove this,  
fix $x\in \R^{N}$. Since $|f(t)|\leq C'(|t|^{p-1}+|t|^{\p-1})$ for all $t\in \R$ (note that $\frac{N+\alpha}{N-p}\frac{p}{2}<\p$ due to $\alpha\in (0, N)$), where $C'>0$ is a constant, $u\in L^{r}(\R^{N})$ for all $r\in [p, \infty]$ and $K=I_{\alpha}*F(u)\in L^{\infty}(\R^{N})$, it follows that 
$|\psi|\leq (1+C'\|K\|_{L^{\infty}(\R^{N})})|u|^{p-1}+f_{*}$ with $f_{*}:=C'\|K\|_{L^{\infty}(\R^{N})} \|u\|^{p-1}_{L^{\infty}(\R^{N})}|u|^{\p-p}\in L^{\frac{N}{p-\e}}(\R^{N})$ for all $\e\in (0, 1]$.
From \cite[Theorem 1]{S}, we derive that there exists a constant $C''>0$, independent of $x$, such that
$$
\|u\|_{L^{\infty}(B_{1}(x))}\leq C''\left(\|u\|_{L^{p}(B_{2}(x))}+\|f_{*}\|^{\frac{1}{p-1}}_{L^{\frac{N}{p-\e}}(B_{2}(x))}\right).
$$
Thus, $|u(x)|\ri 0$ as $|x|\ri \infty$. Since $I_{\alpha}*F(u)\in C_{0}(\R^{N})$ and $f(t)=o(|t|^{p-1})$ as $|t|\ri 0$, we can see that
$(I_{\alpha}*F(u))f(u)\leq \frac{1}{2} u^{p-1}$ in $B_{R}^{c}(0)$, for some $R>0$. Hence, $u$ satisfies $-\Delta_{p}u+\frac{1}{2}u^{p-1}\leq 0$ in $B_{R}^{c}(0)$, and arguing as in the proof of \cite[Lemma 2.3]{AY} we deduce the desired exponential decay of $u$.
\end{remark}
Next we show that $u$ satisfies \eqref{POHOZAEV}. Firstly,  we recall the following useful result.
\begin{lem}\cite{DGMS}\label{DGMSLEMMA}
Let $\Omega\subset \R^{N}$ be an open set, $\mathcal{L}: \Omega\times \R\times \R^{N}\ri \R$ a function of class $C^{1}$ and $g\in L^{\infty}_{loc}(\Omega)$. Assume also that $\xi\mapsto \mathcal{L}(x, s, \xi)$ is strictly convex for each $(x, s)\in \Omega\times \R$. Let $u:\Omega\times \R$ be a locally Lipschitz solution of
\begin{align*}
-\dive\{\nabla_{\xi}\mathcal{L}(x, u, \nabla u)\}+D_{s}\mathcal{L}(x, u, \nabla u)=g \,\, \mbox{ in } \mathcal{D}'(\Omega).
\end{align*}
Then
\begin{align}\begin{split}\label{DGMS}
&\sum_{i, j=1}^{N} \int_{\Omega} D_{i}h_{j}D_{\xi_{i}} \mathcal{L}(x, u, \nabla u) D_{j}u\, dx -\int_{\Omega} \left[ (\dive h) \mathcal{L}(x, u, \nabla u)+h \cdot \nabla_{x} \mathcal{L}(x, u, \nabla u)\right] \, dx \\
&\quad =\int_{\Omega} (h \cdot \nabla u)g\, dx
\end{split}\end{align}
for every $h\in C^{1}_{c}(\Omega; \R^{N})$.
\end{lem}

Take $\mathcal{L}(x, s, \xi)=\frac{1}{p}|\xi|^{p}$ and $g=-|u|^{p-2}u+(I_{\alpha}*F(u))f(u)$ in Lemma \ref{DGMSLEMMA}.
Let $\varphi\in C^{1}_{c}(\R^{N})$ be such that $0\leq \varphi\leq 1$ in $\R^{N}$, $\varphi(x)=1$ for $|x|\leq 1$ and $\varphi(x)=0$ for $|x|\geq 2$. Put
$$
h(x):=\varphi\left(\frac{x}{k}\right)x\in C^{1}(\R^{N}; \R^{N}).
$$
Note that if $h_{j}(x)=\varphi\left(\frac{x}{k}\right)x_{j}$, for $j=1, \dots, N$, then
\begin{align*}
&D_{i}h_{j}(x)=D_{i}\varphi\left(\frac{x}{k}\right) \frac{x_{j}}{k}+\varphi\left(\frac{x}{k} \right)\delta_{ij}, \quad \mbox{ for all } x\in \R^{N}, \, j=1, \dots, N, \\
&\dive h(x)=D\varphi\left(\frac{x}{k}\right)\cdot \frac{x}{k}+N\varphi\left(\frac{x}{k} \right) \quad \mbox{ for all } x\in \R^{N},
\end{align*}
where $\delta_{ij}$ denotes the Kronecker delta symbol. We also observe that
\begin{align}\label{BOUND}
\left|D_{i}\varphi\left(\frac{x}{k}\right) \frac{x_{j}}{k}\right|\leq C \quad \mbox{ for all } x\in \R^{N}, \, i, j=1, \dots, N.
\end{align}
On account of \eqref{DGMS}, we have
\begin{align}\begin{split}\label{DS1}
&\sum_{i, j=1}^{N} \int_{\R^{N}}  D_{i}\varphi\left(\frac{x}{k}\right) \frac{x_{j}}{k} D_{\xi_{i}} \mathcal{L}(x, u, \nabla u) D_{j}u\, dx +\int_{\R^{N}} \varphi\left(\frac{x}{k}\right) D_{\xi}\mathcal{L}(x, u, \nabla u)\cdot \nabla  u\, dx \\
&\quad -\int_{\Omega} \left[ D\varphi\left(\frac{x}{k}\right)\cdot \frac{x}{k} \mathcal{L}(x, u, \nabla u)+N \varphi\left(\frac{x}{k} \right) \mathcal{L}(x, u, \nabla u) \right] \, dx \\
&\quad =\int_{\R^{N}} \left(\varphi\left(\frac{x}{k}\right)x \cdot \nabla u \right)g\, dx.
\end{split}\end{align}
In view of  \eqref{BOUND}, $\varphi\left(\frac{x}{k}\right)\ri 1$ and $D\varphi(\frac{x}{k})\cdot \frac{x}{k}\ri 0$ as $k\ri \infty$, and $|\nabla u|^{p}\in L^{1}(\R^{N})$, it follows from the dominated convergence theorem that
\begin{align}\begin{split}\label{DS2}
&\sum_{i, j=1}^{N} \int_{\R^{N}}  D_{i}\varphi\left(\frac{x}{k}\right) \frac{x_{j}}{k} D_{\xi_{i}} \mathcal{L}(x, u, \nabla u) D_{j}u\, dx +\int_{\R^{N}} \varphi\left(\frac{x}{k}\right) D_{\xi}\mathcal{L}(x, u, \nabla u)\cdot \nabla  u\, dx \\
&\quad -\int_{\Omega} \left[ D\varphi\left(\frac{x}{k}\right)\cdot \frac{x}{k} \mathcal{L}(x, u, \nabla u)+N \varphi\left(\frac{x}{k} \right) \mathcal{L}(x, u, \nabla u) \right] \, dx \\
&\quad\ri \left(1-\frac{N}{p} \right) \int_{\R^{N}} |\nabla u|^{p}\, dx,
\end{split}\end{align}
as $k\ri \infty$.
Now, we note that an integration by parts and the dominated convergence theorem imply that
\begin{align}\begin{split}\label{DS3}
-\int_{\R^{N}} \left(\varphi\left(\frac{x}{k}\right)x \cdot \nabla u \right) |u|^{p-2}u\, dx&=-\int_{\R^{N}} \varphi\left(\frac{x}{k}\right)x \cdot \nabla\left( \frac{|u|^{p}}{p}\right) \, dx \\
&=\int_{\R^{N}} \left[N\varphi\left(\frac{x}{k}\right)+ \nabla \varphi\left( \frac{x}{k}\right)\cdot \frac{x}{k} \right] \frac{|u|^{p}}{p}\, dx\ri \frac{N}{p} \int_{\R^{N}} |u|^{p}\, dx,
\end{split}\end{align}
as $k\ri \infty$. On the other hand, an integration by parts gives
\begin{align*}
&\int_{\R^{N}} \left(\varphi\left(\frac{x}{k}\right)x \cdot \nabla u \right) (I_{\alpha}*F(u)) f(u)\, dx \\
&=\int_{\R^{N}} \int_{\R^{N}} (F\circ u)(y) I_{\alpha}(x-y) \varphi\left(\frac{x}{k} \right)x \cdot \nabla (F\circ u)(x)\, dx dy \\
&=\frac{1}{2} \int_{\R^{N}} \int_{\R^{N}} I_{\alpha}(x-y) \left( (F\circ u)(y) \varphi\left(\frac{x}{k} \right)x\cdot \nabla (F\circ u)(x)+(F\circ u)(x) \varphi\left(\frac{y}{k} \right)y\cdot \nabla (F\circ u)(y)  \right)\, dx dy \\
&=-\int_{\R^{N}} \int_{\R^{N}} F(u(y)) I_{\alpha}(x-y) \left[N\varphi\left( \frac{x}{k}\right)+x\cdot \nabla \varphi\left( \frac{x}{k} \right)  \right] F(u(x))\, dx dy\\
&\quad+\frac{N-\alpha}{2} \int_{\R^{N}} \int_{\R^{N}} F(u(y)) I_{\alpha}(x-y) \frac{(x-y)\cdot (x\varphi\left(\frac{x}{k}\right)- y\varphi\left(\frac{y}{k}\right) )}{|x-y|^{2}} F(u(x))\, dx dy,
\end{align*}
and applying the dominated convergence theorem we see that 
\begin{align}\label{DS4}
\int_{\R^{N}} \left(\varphi\left(\frac{x}{k}\right)x \cdot \nabla u \right) (I_{\alpha}*F(u)) f(u)\, dx\ri -\frac{N+\alpha}{2} \int_{\R^{N}} (I_{\alpha}*F(u))F(u)\, dx,
\end{align}
as $k\ri \infty$. Combining \eqref{DS1}, \eqref{DS2}, \eqref{DS3} and \eqref{DS4}, we obtain \eqref{POHOZAEV}.
The proof of Theorem \ref{thm1} is now completed.

\end{document}